\documentclass{amsart}
\usepackage{amssymb}
\usepackage[foot]{amsaddr}
\usepackage{mathtools}
\usepackage[margin=1in]{geometry}
\usepackage{biblatex}
\usepackage{todonotes}
\usepackage{enumitem}
\usepackage[colorlinks=true]{hyperref}

\renewcommand{\v}{\mathbf{v}}
\newcommand{\w}{\mathbf{w}}
\newcommand{\x}{\mathbf{x}}
\newcommand{\y}{\mathbf{y}}
\newcommand{\z}{\mathbf{z}}

\newtheorem{theorem}{Theorem}
\newtheorem{lemma}{Lemma}
\newtheorem{proposition}{Proposition}
\theoremstyle{remark}
\newtheorem{remark}{Remark}

\title[Convergence of Extended Douglas--Rachford Splitting]{Sharp Convergence of Extended Douglas--Rachford Splitting with Variable Relaxation}
\author{Matteo Santamaria}
\address{Industrial Engineering \& Operations Research \\ University of California, Berkeley}
\date{November 2026}

\begin{document}

\begin{abstract}
    We establish a sharp convergence theory for extended Douglas--Rachford splitting when the first operator is maximally monotone and the second is a convex subdifferential on a real Hilbert space. The resolvent stepsizes may be chosen independently, and the relaxation parameter may depend on the current and past iterates. For fixed stepsizes $\alpha,\beta>0$, let $\vartheta=\min\{2,2\alpha/\beta\}$. Assuming the sum of the operators has a zero, we prove weak convergence under $0<\theta_k\le\vartheta$ and $\sum_k\theta_k(\vartheta-\theta_k)=\infty$. Within these pointwise bounds, the divergence condition is also necessary for a universal convergence guarantee when the relaxation sequence is prescribed in advance. Constant relaxation recovers the exact range $0<\theta<\vartheta$, resolving a conjecture of Nilsson, Åkerman, and Giselsson. We also establish a sharp bound on growth of the residual, defined as the difference between the two resolvent outputs. We combine this bound with a new Lyapunov function whose definition is independent of the relaxation parameter to obtain a last-iterate residual rate, including $o(N^{-1/2})$ when the relaxation parameters remain uniformly separated from both endpoints.
\end{abstract}

\maketitle

\vspace{-18pt}

\section{Introduction}

We study the monotone inclusion
\begin{equation}
    \label{eq:monotone_inclusion}
    \text{find} \quad \x \in \mathcal{H} \quad \text{such that} \quad 0 \in F(\x) + \partial g(\x),
\end{equation}
where $\mathcal{H}$ is a nonzero real Hilbert space, $F:\mathcal{H}\rightrightarrows\mathcal{H}$ is maximally monotone, and $g:\mathcal{H}\to(-\infty,+\infty]$ is proper, lower semicontinuous, and convex. Many problems in optimization and economics admit a formulation of this kind, including convex minimization and the computation of Nash equilibria in monotone games \cite{briceno2013nash}. The following examples illustrate two important instances.

\begin{enumerate}[label=\roman*.]
    \item \textbf{Convex optimization.} Let $f:\mathcal{H}\to(-\infty,+\infty]$ be proper, lower semicontinuous, and convex. Consider
    \[
        \min_{\x\in\mathcal{H}}\; f(\x)+g(\x).
    \]
    With \(F=\partial f\), the inclusion becomes
    \[
    0\in\partial f(\x)+\partial g(\x),
    \]
    which is sufficient for \(\x\) to minimize \(f+g\), and necessary under a suitable subdifferential sum qualification \cite[Theorem~16.3 and Corollary~16.48]{bauschke2017convex}.

    \item \textbf{Constrained variational inequalities.} Suppose $F:\mathcal{H}\to\mathcal{H}$ is single-valued, and let $C\subseteq\mathcal{H}$ be nonempty, closed, and convex. Consider the variational inequality
    \[
        \text{find} \quad \x \in C \quad \text{such that} \quad \langle F(\x) , u - \x \rangle \ge 0 \quad \text{for every} \ u \in C.
    \]
    With $g=\iota_C$, the indicator of $C$, this problem is equivalent to the inclusion~\eqref{eq:monotone_inclusion}.
    Nash equilibrium conditions in monotone games provide one application of this variational inequality.
    Here, $F$ need not be the subdifferential of a convex function,
    so the formulation extends beyond convex minimization.
\end{enumerate}

To solve \eqref{eq:monotone_inclusion}, we consider a version of extended Douglas--Rachford splitting \cite{nilsson2025extendingdouglasrachfordsplittingconvex} with adaptive relaxation. Each iteration treats the two operators separately, evaluating the resolvent of $F$ followed by the proximal map of $g$. For fixed resolvent stepsizes $\alpha,\beta>0$ and an initial point $\z^0\in\mathcal{H}$, the iteration is
\begin{subequations}
\label{eq:extended_douglas_rachford}
\begin{align}
    \x^k &\gets J_{\alpha F}(\z^k) \\
    \y^k &\gets \operatorname{prox}_{\beta g} \left( \x^k + \frac{\beta}{\alpha} \left( \x^k - \z^k \right) \right) \\
    \z^{k + 1} &\gets \z^k + \theta_k \left( \y^k - \x^k \right).
\end{align}
\end{subequations}

We call the difference $\y^k-\x^k$ between the two resolvent outputs the \emph{residual} at iteration $k$. When it vanishes, their common value solves the inclusion. Since both outputs are available before the update of \(\z^{k+1}\), the relaxation parameter \(\theta_k\) may be prescribed in advance or chosen using the current and past residuals. 

When $\alpha=\beta$, the iteration reduces to a relaxed form of the Douglas--Rachford method originally studied by Lions and Mercier \cite{lions1979splitting}. In this setting, Combettes proved that, whenever the inclusion has a solution, \((\z^k)\) converges weakly if \(0<\theta_k<2\) and \(\sum_k\theta_k(2-\theta_k)=\infty\) \cite[Corollary~5.2]{combettes2004solving}. See also \cite[Theorem~26.11]{bauschke2017convex}.

In dimensions at least two, no fixed relaxation parameter guarantees convergence for every pair of maximally monotone operators when the resolvent stepsizes are unequal \cite[Theorem~2]{ryu2020uniqueness}. Nilsson, Åkerman, and Giselsson established a sharp range of constant relaxation parameters when the first operator is a convex subdifferential and the second is any maximally monotone operator. Whether the same bounds suffice when only the second operator is assumed to be a convex subdifferential remained open \cite[Theorem~3.3 and Conjecture~7.1]{nilsson2025extendingdouglasrachfordsplittingconvex}.

Our convergence theorem permits variable relaxation and resolves this ordering question as a constant-relaxation specialization. Whenever the inclusion has a solution, our conditions guarantee weak convergence with independently chosen resolvent stepsizes and allow the relaxation parameter to depend on the current and past iterates.

For the variational inequality above, our theorem justifies
evaluating the resolvent of $F$ before projecting onto $C$.
The result of Nilsson et al. covers the opposite ordering,
with the projection first. Exchanging the operators leaves
the inclusion unchanged but changes the iteration, so
convergence for one ordering does not directly establish
convergence for the other.

\subsection{Contributions}

Our contributions comprise a convergence theorem for variable relaxation, its consequences for constant relaxation and sharpness, and a sharp residual stability bound with a last-iterate rate.

\begin{enumerate}[label=\roman*.]
    \item \textbf{Weak convergence with variable relaxation.} For fixed, possibly unequal resolvent stepsizes, we give conditions on the relaxation sequence under which $(\z^k)$ converges weakly to a fixed point and both resolvent outputs converge weakly to the same solution whenever the inclusion has a solution (Theorem~\ref{thm:edr_convergence}). The relaxation parameter may depend on the iterates. Our divergence condition extends the classical one for relaxed Douglas--Rachford splitting by replacing $2$ with $\min\{2,2\alpha/\beta\}$, which recovers the condition of Combettes when the resolvent stepsizes coincide \cite[Corollary~5.2]{combettes2004solving}.

    \item \textbf{Resolution of the constant-relaxation conjecture.} For constant relaxation, the theorem yields the exact range $0<\theta<\min\{2,2\alpha/\beta\}$, resolving a conjecture of Nilsson, Åkerman, and Giselsson \cite[Conjecture~7.1]{nilsson2025extendingdouglasrachfordsplittingconvex}. For relaxation sequences prescribed in advance within the theorem's pointwise bounds, the divergence condition is also necessary: if it fails, we construct a problem satisfying our assumptions with a unique solution for which $(\x^k)$ does not converge to that solution.

    \item \textbf{Sharp residual stability and a last-iterate rate.} Between any two iterations, the squared residual can grow by at
    most the optimal factor $\max\{\alpha/\beta,\beta/\alpha\}$
    (Lemma~\ref{lem:edr_residual_stability}). This factor depends
    only on the stepsize ratio and does not accumulate with the
    number of intervening iterations. With equal stepsizes, the residual norm is therefore nonincreasing. The bound holds under the theorem's pointwise relaxation bounds without assuming that the inclusion has a solution or that the divergence condition holds.

    Under the full hypotheses of Theorem~\ref{thm:edr_convergence}, it yields a last-iterate rate that depends on the relaxation sequence. In particular, the last-iterate residual satisfies $\|\y^N-\x^N\|=o(N^{-1/2})$ whenever the relaxation parameters stay a fixed positive distance from both endpoints of the admissible interval, even if they are chosen adaptively.

\end{enumerate}

Unequal resolvent stepsizes can destroy averagedness of the update \cite[discussion following Theorem~3.3]{nilsson2025extendingdouglasrachfordsplittingconvex} and monotonicity of the residual norm, so the classical convergence and last-iterate arguments do not apply directly. Our analysis uses a Lyapunov function whose definition is independent of the relaxation parameter to obtain descent as the relaxation varies. A second, nonincreasing functional is uniformly comparable to the squared residual and controls later residuals in terms of earlier ones. The uniform residual bound converts the weighted estimate supplied by Lyapunov descent into a last-iterate rate.

\subsection{Related literature}

Ryu characterized fixed-parameter resolvent splittings for sums of two maximally monotone operators on $\mathbb{R}^d$. He considered methods with a single state vector in $\mathbb{R}^d$ that evaluate each resolvent once per iteration and otherwise use linear combinations with fixed scalar coefficients. Requiring a solution map to send the fixed-point set onto the solution set for every operator pair, he showed that all such splittings have the update in \eqref{eq:extended_douglas_rachford}, up to rescaling the state and exchanging the operators, with fixed $\alpha,\beta>0$ and constant $\theta\ne0$ \cite[Theorem~1]{ryu2020uniqueness}.

Convergence imposes further restrictions. For $d\ge2$, Ryu proved that these splittings converge from every initial point for every maximally monotone pair whose sum has a zero if and only if $\alpha=\beta$ and $0<\theta<2$ \cite[Theorem~2]{ryu2020uniqueness}. This leaves open which additional properties of the operators permit unequal resolvent stepsizes.

Dao and Phan obtained such a guarantee through an adaptive Douglas--Rachford method whose iteration is equivalent, after reparameterization, to the constant-relaxation update studied here \cite[equations~(7)--(8)]{nilsson2025extendingdouglasrachfordsplittingconvex}. They chose fixed parameters according to the operators' monotonicity constants and proved weak convergence of $(\z^k)$ to a fixed point, allowing unequal stepsizes when, for example, strong monotonicity of one operator compensates for weak monotonicity of the other \cite[condition~(86) and Theorem~4.5]{dao2019adaptive}. Here, ``adaptive'' refers to the dependence of the parameters on the operator constants.

Bartz, Dao, and Phan used conical averagedness to enlarge the range of fixed parameters guaranteeing weak convergence of $(\z^k)$ when the monotonicity constants have positive sum \cite[Theorem~5.7]{bartz2022conical}. Complementing this parameter analysis, Bartz, Campoy, and Phan established weak convergence of the first resolvent output $\x^k=J_{\alpha F}(\z^k)$ to a solution under the corresponding parameter conditions, including when the monotonicity constants sum to zero \cite[Theorem~4.1]{bartz2020demiclosedness}.

When both monotonicity constants are zero, however, these convergence conditions still require equal resolvent stepsizes, even for two convex subdifferentials. This left open whether convex subdifferential structure alone could permit unequal resolvent stepsizes.

Nilsson, Åkerman, and Giselsson showed that convex subdifferential structure permits unequal resolvent stepsizes without strong monotonicity. For two convex subdifferentials and arbitrary fixed stepsizes $\alpha,\beta>0$, they identified the exact range of constant relaxation parameters guaranteeing weak convergence whenever the inclusion has a solution:
$0<\theta<\min\{2,2\alpha/\beta\}$ \cite[Theorem~3.3]{nilsson2025extendingdouglasrachfordsplittingconvex}. Their Lyapunov analysis in fact requires subdifferential structure only for the operator evaluated first; the second operator may be any maximally monotone operator.

The reverse ordering remained open when only the second operator is assumed to be a convex subdifferential. Their counterexamples with two subdifferentials already establish necessity of the same parameter bounds for this larger class. Conjecture~7.1 asserts that these bounds are also sufficient \cite[Section~4.2 and Conjecture~7.1]{nilsson2025extendingdouglasrachfordsplittingconvex}.

Their Lyapunov function depends explicitly on the fixed relaxation
parameter \cite[Definition~4.2]{nilsson2025extendingdouglasrachfordsplittingconvex}.
Allowing the relaxation to vary changes the Lyapunov function
between iterations, so their fixed-parameter descent identity
does not directly yield a telescoping estimate.

Adaptive parameter choices in Douglas--Rachford splitting include changes to the resolvent stepsize itself. Lorenz and Tran-Dinh vary the common resolvent stepsize across iterations and prove convergence when it converges to a positive limit and its successive changes are absolutely summable \cite{lorenz2019nonstationary}. More recently, Atenas, Bauschke, Dao, and Tam introduced a fixed-point relocation step to account for changes in the fixed-point set as the resolvent stepsize varies \cite{atenas2026relocated}. Their framework yields convergent variable-stepsize Douglas--Rachford and multioperator splitting methods.

For relaxation of the fixed-point update, Giselsson, Fält, and Boyd develop a line search that tests the residual at trial points along the update direction. This requires trial operator evaluations, although their additional cost can be small for particular operator structures \cite{giselsson2016linesearch}. Iutzeler and Hendrickx develop online relaxation rules within an averaged-operator framework \cite{iutzeler2019online}. For Douglas--Rachford splitting with unequal resolvent stepsizes, Briceño-Arias, Roldán, and Taylor establish averagedness and convergence with variable relaxation under cocoercivity assumptions \cite{briceno2026relaxed}.

Our theorem requires no cocoercivity assumption, and the relaxation parameters need not converge or have summable successive changes.

Last-iterate residual rates are established for classical Douglas--Rachford splitting. Davis and Yin derived residual bounds for averaged iterations, yielding $\|\y^N-\x^N\|=o(N^{-1/2})$ for classical Douglas--Rachford splitting when the relaxation parameters stay uniformly away from the endpoints of the admissible interval \cite{davis2016convergence}. Dao and Phan obtained the same residual rate for their fixed-parameter adaptive Douglas--Rachford method under the monotonicity and parameter conditions discussed above \cite[Theorem~4.5]{dao2019adaptive}.

With unequal resolvent stepsizes, the residual norm need not decrease. Our sharp stability estimate (Lemma~\ref{lem:edr_residual_stability}) replaces the residual monotonicity used in the classical last-iterate analysis.

\section{Background and Preliminaries}

We denote the inner product on $\mathcal{H}$ by $\langle \cdot, \cdot \rangle$ and its induced norm by $\| \cdot \|$. We briefly review the concepts and standard results used in the convergence analysis. Unless otherwise indicated, we refer to \cite{bauschke2017convex} for proofs and further background.

\subsection{Monotone operator theory}

For a set-valued operator $F : \mathcal{H}\rightrightarrows\mathcal{H}$, its inverse relation is defined by
\[
    F^{-1}(w)\coloneq\{x\in\mathcal{H}:w\in F(x)\}.
\]
The zero set of $F$ is
\[
    \operatorname{zer}F\coloneq\{x\in\mathcal{H}:0\in F(x)\}.
\]

An operator $F$ is monotone if
\[
    \langle x_1-x_2,w_1-w_2\rangle\ge0
    \qquad\text{whenever }w_i\in F(x_i),\quad i=1,2,
\]
and maximally monotone if its graph has no proper monotone extension. For maximally monotone $F$, membership in its graph is characterized by
\begin{equation}
    \label{eq:maximal_graph_characterization}
    u\in F(x)\quad\Longleftrightarrow\quad
    \langle x-y,u-v\rangle\ge0
    \quad\text{for every }y\in\mathcal{H}\text{ and }v\in F(y).
\end{equation}
A pair $(x,u)$ satisfying the inequality on the right is said
to be monotonically related to the graph of $F$.
Thus maximal monotonicity ensures that every such pair
already belongs to the graph.

Inversion preserves monotonicity and maximal monotonicity.
If $F,G:\mathcal{H}\rightrightarrows\mathcal{H}$ are maximally monotone, then their Cartesian product, defined by
\[
    (F\times G)(x,y)\coloneq F(x)\times G(y),
\]
is maximally monotone on the product Hilbert space $\mathcal{H}\times\mathcal{H}$ \cite[Proposition~23.18]{bauschke2017convex}.

For a monotone operator $F$ and fixed $\lambda>0$, Minty's theorem \cite[Theorem~21.1]{bauschke2017convex} gives
\[
    F\text{ is maximally monotone}
    \quad\Longleftrightarrow\quad
    \operatorname{ran}(I+\lambda F)=\mathcal{H}.
\]
If $F$ is maximally monotone and $\lambda>0$, the resolvent $J_{\lambda F}\coloneq(I+\lambda F)^{-1}$ is a single-valued mapping from $\mathcal{H}$ to $\mathcal{H}$, characterized by
\begin{equation}
    \label{eq:resolvent_characterization}
    x=J_{\lambda F}(z)\quad\Longleftrightarrow\quad z-x\in\lambda F(x).
\end{equation}
The resolvent $J_{\lambda F}$ is firmly nonexpansive, meaning that
\[
    \|J_{\lambda F}(u)-J_{\lambda F}(v)\|^2\le\langle J_{\lambda F}(u)-J_{\lambda F}(v),u-v\rangle
    \qquad\text{for all }u,v\in\mathcal{H}.
\]
Consequently, $J_{\lambda F}$ is nonexpansive.

\subsection{Convex functions and proximal operators}

Let $g:\mathcal{H}\to(-\infty,+\infty]$ be proper, lower semicontinuous, and convex. Its effective domain is
\[
    \operatorname{dom}g\coloneq\{x\in\mathcal{H}:g(x)<+\infty\}.
\]
For $x\in\operatorname{dom}g$, its subdifferential is
\[
    \partial g(x)\coloneq\{v\in\mathcal{H}:g(y)\ge g(x)+\langle v,y-x\rangle\text{ for all }y\in\mathcal{H}\},
\]
and $\partial g(x)=\varnothing$ when $x\notin\operatorname{dom}g$. The operator $\partial g$ is maximally monotone.

For $\lambda>0$, the proximal operator is
\[
    \operatorname{prox}_{\lambda g}(z)\coloneq\underset{x\in\mathcal{H}}{\operatorname{argmin}}\left\{g(x)+\frac{1}{2\lambda}\|x-z\|^2\right\}=J_{\lambda\partial g}(z).
\]
As a resolvent, it is single-valued, defined on all of $\mathcal{H}$, and nonexpansive. Its optimality condition is
\begin{equation}
    \label{eq:proximal_optimality}
    x=\operatorname{prox}_{\lambda g}(z)\quad\Longleftrightarrow\quad\frac{z-x}{\lambda}\in\partial g(x).
\end{equation}

For a nonempty closed convex set $C\subseteq\mathcal{H}$, its indicator function $\iota_C$ is zero on $C$ and $+\infty$ elsewhere. For every $\lambda>0$, $\operatorname{prox}_{\lambda\iota_C}=P_C$, where $P_C$ denotes the orthogonal projection onto $C$.

\subsection{Generalized Bregman divergences}

Let $\phi:\mathcal{H}\to(-\infty,+\infty]$ be proper, lower semicontinuous, and convex. Fix a base point $y\in\mathcal{H}$ and a subgradient $v\in\partial\phi(y)$. The generalized Bregman divergence from $y$ to $x$, relative to $v$, is \cite[Definition~2.2]{burger2015bregmandistancesinverseproblems}
\begin{equation*}
    D_\phi(x,y,v)\coloneq\phi(x)-\phi(y)-\langle v,x-y\rangle.
\end{equation*}
This generalized definition requires neither the differentiability nor the strict convexity commonly assumed in the classical definition. By the subgradient inequality, $D_\phi(x,y,v)\ge0$ for all $x\in\mathcal{H}$, and it is finite exactly when $x\in\operatorname{dom}\phi$. For fixed $y$ and $v$, the map $x\mapsto D_\phi(x,y,v)$ is convex, since it differs from $\phi$ by an affine function.

The divergence is invariant under translation and addition of an affine function. Let $a,b\in\mathcal{H}$ and $c\in\mathbb{R}$, and define
\[
    h(x)\coloneq\phi(a+x)+\langle b,x\rangle+c.
\]
The subgradient inequality gives $\partial h(y)=\partial\phi(a+y)+b$. Hence, for $v\in\partial\phi(a+y)$ and $a+x\in\operatorname{dom}\phi$, the affine terms cancel in the divergence, giving
\begin{equation}
    \label{eq:bregman_affine_translation}
    D_h(x,y,v+b)=D_\phi(a+x,a+y,v).
\end{equation}

\subsection{Fenchel conjugacy}
\label{sec:fenchel_conjugacy}

For a proper, lower semicontinuous, convex function $g:\mathcal{H}\to(-\infty,+\infty]$, its Fenchel conjugate is
\[
    g^*(v)\coloneq\sup_{x\in\mathcal{H}}\{\langle v,x\rangle-g(x)\}.
\]
Both $g^*$ and $u\mapsto g^*(-u)$ are proper, lower semicontinuous, and convex. The Fenchel--Young inequality states that
\[
    g(x)+g^*(v)\ge\langle v,x\rangle
\]
for all $x,v\in\mathcal{H}$. Equality is characterized by \cite[Theorem~16.29]{bauschke2017convex}
\begin{equation}
    \label{eq:fenchel_subgradient}
    v\in\partial g(x)\quad\Longleftrightarrow\quad g(x)+g^*(v)=\langle v,x\rangle\quad\Longleftrightarrow\quad x\in\partial g^*(v).
\end{equation}

\subsection{Weak convergence}
\label{sec:weak_convergence}

We write $x^k\rightharpoonup x$ if $\langle x^k,u\rangle\to\langle x,u\rangle$ for every $u\in\mathcal{H}$. A point $u\in\mathcal{H}$ is a weak cluster point of $(u^k)$ if $u^{k_j}\rightharpoonup u$ along some subsequence. Every bounded sequence in $\mathcal{H}$ has a weakly convergent subsequence. Hence, if $u$ is the unique weak cluster point of a bounded sequence $(u^k)$, then $u^k\rightharpoonup u$.

\section{Convergence Results}

The following theorem gives conditions on the relaxation sequence in \eqref{eq:extended_douglas_rachford} that ensure weak convergence and quantify the decay of the last-iterate residual.
The corresponding necessity results are established in
Proposition~\ref{prop:edr_parameter_sharpness}.

\begin{theorem} \label{thm:edr_convergence}
    Suppose the inclusion~\eqref{eq:monotone_inclusion} has a solution. Fix $\alpha,\beta>0$. If the relaxation sequence $(\theta_k)_{k\ge0}$ satisfies
    \begin{equation}
        \label{eq:edr-parameters}
        0<\theta_k\le\vartheta\coloneq\min\left\{2,\frac{2\alpha}{\beta}\right\}\quad(k\ge0),\qquad
        \sum_{k=0}^{\infty}\theta_k(\vartheta-\theta_k)=\infty,
    \end{equation}
    then for every $\z^0 \in \mathcal{H}$, there are $\x^\infty, \z^\infty$ such that
    \[
        (\x^k,\y^k,\z^k)
        \rightharpoonup(\x^\infty,\x^\infty,\z^\infty),
    \]
    with
    \begin{equation}
        \label{eq:edr_limit_identities}
        \x^\infty=J_{\alpha F}(\z^\infty)
        =\operatorname{prox}_{\beta g}
        \left(\x^\infty+\frac{\beta}{\alpha}
        (\x^\infty-\z^\infty)\right).
    \end{equation}
    In particular, $0\in F(\x^\infty)+\partial g(\x^\infty)$, and the last-iterate residual satisfies
    \begin{equation}
        \label{eq:edr_residual_rate}
        \left(\sum_{k=0}^{N-1}\theta_k(\vartheta-\theta_k)\right)\|\y^N-\x^N\|^2\to0.
    \end{equation}
\end{theorem}

We first record three immediate consequences of Theorem~\ref{thm:edr_convergence}.

\begin{remark}[Residual rate]
    Under the hypotheses of Theorem~\ref{thm:edr_convergence},
    if the relaxation parameters stay uniformly away from $0$
    and $\vartheta$, the sum in \eqref{eq:edr_residual_rate}
    grows linearly with $N$, giving
    $\|\y^N-\x^N\|=o(N^{-1/2})$.
    This includes every constant relaxation
    $0<\theta<\vartheta$.
\end{remark}

\begin{remark}[History-dependent relaxation]
    \label{rem:adaptive_relaxation}
    Suppose \eqref{eq:monotone_inclusion} has a solution.
    Prescribe $\theta_k=\vartheta/2$ at iterations
    $k=1,2,4,8,\ldots$, and at every remaining iteration
    choose any $\theta_k\in(0,\vartheta]$ using the current
    and past iterates.

    Each prescribed midpoint step contributes $\vartheta^2/4$
    to the damping sum, and all other contributions are nonnegative.
    Hence, for $N\ge2$,
    \[
        \sum_{k=0}^{N-1}\theta_k(\vartheta-\theta_k)
        \ge\frac{\vartheta^2}{4}\lfloor\log_2 N\rfloor
        \to\infty.
    \]
    Theorem~\ref{thm:edr_convergence} therefore applies.
    Among the first $N$ iterations, only $O(\log N)$ relaxation
    parameters are prescribed, so their proportion tends to zero.
\end{remark}

\begin{remark}[Fixed metric]
    \label{rem:edr_fixed_metric}
    The convergence theorem also permits a fixed metric in both resolvent evaluations, with the same relaxation conditions.
    Let $M:\mathcal{H}\to\mathcal{H}$ be a bounded self-adjoint linear operator satisfying $\langle Mx,x\rangle\ge c\|x\|^2$ for some $c>0$ and all $x\in\mathcal{H}$. Write $\mathcal{H}_M$ for $\mathcal{H}$ equipped with the inner product $\langle x,y\rangle_M\coloneq\langle Mx,y\rangle$ and the induced norm $\|x\|_M\coloneq\sqrt{\langle Mx,x\rangle}$. For $\lambda>0$, define the metric resolvent and proximal map by
    \[
    \begin{aligned}
        J^M_{\lambda F} &\coloneq (I+\lambda M^{-1}F)^{-1},\\
        \operatorname{prox}^M_{\lambda g}(z)
        &\coloneq \underset{x\in\mathcal{H}}{\operatorname{argmin}}
        \left\{g(x)+\frac{1}{2\lambda}\|x-z\|_M^2\right\}
        =J^M_{\lambda\partial g}(z).
    \end{aligned}
    \]
    For fixed stepsizes $\alpha,\beta>0$ and an initial point $\z^0\in\mathcal{H}$, the fixed-metric iteration is
    \[
    \begin{aligned}
        \x^k &\gets J^M_{\alpha F}(\z^k),\\
        \y^k &\gets \operatorname{prox}^M_{\beta g}
            \left(\x^k+\frac{\beta}{\alpha}(\x^k-\z^k)\right),\\
        \z^{k+1} &\gets \z^k+\theta_k(\y^k-\x^k).
    \end{aligned}
    \]

    Under the hypotheses of Theorem~\ref{thm:edr_convergence}, the weak convergence conclusions hold for this iteration, and the last-iterate residual rate holds in $\|\cdot\|_M$. These conclusions follow by applying Theorem~\ref{thm:edr_convergence} in $\mathcal{H}_M$, where $M^{-1}F$ is maximally monotone and the subdifferential of $g$ is $M^{-1}\partial g$. The inclusion has the same solutions, and equivalence of the two norms gives the same weak topology.
\end{remark}

\subsection{Analysis of the normalized iteration}

To prove Theorem~\ref{thm:edr_convergence}, we return to \eqref{eq:extended_douglas_rachford} and first analyze it with stepsizes $(\alpha,\beta)=(1,\gamma)$, where $\gamma>0$. This case suffices for arbitrary stepsizes $\alpha,\beta>0$, since applying the normalized iteration to $(\alpha F,\alpha g)$ with $\gamma=\beta/\alpha$ gives the same iterates from the same initial point and relaxation sequence. For a general pair $(F,g)$, the normalized iteration reads
\begin{subequations}
\label{eq:edr_alpha_one}
\begin{align}
    \x^k &\gets J_{F}(\z^k) \label{eq:jf} \\
    \y^k &\gets \operatorname{prox}_{\gamma g} \left( \x^k + \gamma \left( \x^k - \z^k \right) \right) \\
    \z^{k + 1} &\gets \z^k + \theta_k \left( \y^k - \x^k \right). \label{alg:step_3}
\end{align}
\end{subequations}
The resolvent and proximal updates give the corresponding operator values
\begin{equation}
    \label{eq:edr_operator_values}
    \w^k \coloneq \z^k - \x^k \in F(\x^k), \qquad \v^k \coloneq \frac{\x^k - \y^k}{\gamma} - \w^k \in \partial g(\y^k).
\end{equation}

For $0<\gamma\le1$, Lemma~\ref{lem:edr_descent} gives a descent identity for a Lyapunov function whose definition is independent of the relaxation parameter. This allows us to compare successive iterations when $\theta_k$ varies and bound a weighted sum of squared residuals. Lemma~\ref{lem:edr_duality} reduces the case $\gamma>1$ to this case through a primal--dual transformation. Finally, Lemma~\ref{lem:edr_residual_stability} controls subsequent residuals in terms of earlier ones, converting the weighted bound into a last-iterate rate.

\begin{lemma}[Lyapunov descent]
     \label{lem:edr_descent}
     Fix a solution pair satisfying $\w^\star\in F(\x^\star)$ and $-\w^\star\in\partial g(\x^\star)$. Fix $0 < \gamma \le 1$, let $(\theta_k)_{k\ge0}$ satisfy $0<\theta_k\le2$ for every $k\ge0$, and let the sequence $(\x^k, \y^k, \z^k)_{k \ge 0}$ be generated by \eqref{eq:edr_alpha_one}. Then the Lyapunov function
     \begin{equation}
     \begin{aligned}
        \label{eq:lyapunov}
        V_k^\star &\coloneq \frac{1}{\gamma}  \| \x^k - \x^\star \|^2 + \| \w^k - \w^\star + \left( \frac{1}{\gamma} - 1 \right) \left( \y^k - \x^k \right) \|^2 \\ &\qquad + \frac{2}{\gamma} \langle \x^k - \x^\star, \w^k - \w^\star \rangle + \left( \frac{1}{\gamma} - 1 \right) \| \y^k - \x^k \|^2 + 2 \left( \frac{1}{\gamma} - 1 \right) D_{g} \! \left(\y^k, \x^\star, -\w^\star\right).
     \end{aligned}
     \end{equation}
     is nonnegative for all $k \ge 0$, and the sequence $(V_{k}^\star)_{k \ge 0}$
     satisfies the exact descent identity
     \begin{equation}
        \label{eq:exact_descent_identity}
        \begin{aligned}
            V_k^\star - V_{k + 1}^\star &= \theta_k ( 2 - \theta_k) \| \y^k - \x^k \|^2 + \frac{(2 - \theta_k)(1 - \gamma^2)}{\theta_k \gamma^2} \| \x^{k + 1} - \x^k \|^2 \\ &\qquad + \frac{1 - \gamma}{\gamma^2} \| \y^{k + 1} -\y^k \|^2 + 2\theta_k \langle \w^k - \w^\star, \x^k - \x^\star \rangle + 2 \theta_k \langle \v^k + \w^\star,  \y^k - \x^\star \rangle \\ &\qquad + \frac{2(1 - \gamma^2)}{\theta_k \gamma^2} \langle \w^{k + 1} - \w^k, \x^{k + 1} - \x^k \rangle
         + 2 \left( \frac{1}{\gamma} - 1 \right) D_{g} \! \left(\y^k, \y^{k + 1}, \v^{k + 1}\right).
        \end{aligned}
     \end{equation}
     In particular, for every $k\ge0$,
     \begin{equation}
         \label{eq:descent_bound}
         V_k^\star-V_{k+1}^\star\ge\theta_k(2-\theta_k)\|\y^k-\x^k\|^2 \ge0.
     \end{equation}
\end{lemma}
\begin{proof}
    Monotonicity of $F$ implies $\langle \x^k - \x^\star, \w^k - \w^\star \rangle \ge 0$ and $\langle \x^{k + 1} - \x^k, \w^{k + 1} - \w^k \rangle \ge 0$, while monotonicity of $\partial g$ implies $\langle \v^k + \w^\star, \y^k - \x^\star \rangle \ge 0$. Together with the nonnegativity of the generalized Bregman divergence and the coefficients in \eqref{eq:lyapunov} for $0<\gamma\le1$, these inequalities give $V_k^\star\ge0$.

    To derive the descent identity, we translate the chosen solution pair $(\x^\star,\w^\star)$ to the origin. This simplifies the solution-dependent terms in the Lyapunov function and expresses its Bregman term as a value of the shifted function. Define
    \[\bar{F}(u) \coloneq F(\x^\star + u) - \w^\star, \qquad \bar{g}(u) \coloneq g(\x^\star + u) - g(\x^\star) + \langle \w^\star, u \rangle.\]
    Then $0 \in \bar{F}(0), 0 \in \partial \bar{g}(0)$, and $\bar{g}(0) = 0$. Apply the iteration \eqref{eq:edr_alpha_one} to $\bar{F}, \bar{g}$ with the same relaxation sequence $(\theta_k)_{k\ge0}$, starting from $\bar{\z}^0 \coloneq \z^0 - \x^\star - \w^\star$. The identities
    \begin{align*}
        J_{\bar{F}}(u) &= J_{F}(u + \x^\star + \w^\star) - \x^\star,\\
        \operatorname{prox}_{\gamma\bar{g}}(u)
        &= \operatorname{prox}_{\gamma g}(u + \x^\star - \gamma\w^\star) - \x^\star
    \end{align*}
    follow from the resolvent characterization and the proximal optimality condition. Substituting these identities into the iteration and using induction gives, for every $k\ge0$,
    \[
        \bar{\x}^k=\x^k-\x^\star,\qquad
        \bar{\y}^k=\y^k-\x^\star,\qquad
        \bar{\z}^k=\z^k-\x^\star-\w^\star.
    \]

    Define the corresponding operator values by
    \begin{align*}
        \bar{\w}^k &\coloneq \bar{\z}^k - \bar{\x}^k
            = \w^k - \w^\star \in \bar{F}(\bar{\x}^k),\\
        \bar{\v}^k &\coloneq \frac{\bar{\x}^k - \bar{\y}^k}{\gamma} - \bar{\w}^k
            = \v^k + \w^\star \in \partial\bar{g}(\bar{\y}^k).
    \end{align*}
    The Bregman term in the Lyapunov function satisfies
    \[
        D_{g}(\y^k,\x^\star,-\w^\star)
        = g(\y^k) - g(\x^\star) + \langle \w^\star,\y^k-\x^\star\rangle
        = \bar{g}(\bar{\y}^k).
    \]
    Using $\bar{\z}^k=\bar{\x}^k+\bar{\w}^k$ and
    $\bar{\w}^k-(\bar{\y}^k-\bar{\x}^k)=\bar{\z}^k-\bar{\y}^k$,
    expanding the squared norms in \eqref{eq:lyapunov} gives
    \begin{align*}
        V_k^\star
        &= \|\bar{\z}^k\|^2
        +2\left(\frac{1}{\gamma}-1\right)
            \langle\bar{\z}^k,\bar{\y}^k\rangle
        -\left(\frac{1}{\gamma}-1\right)\|\bar{\y}^k\|^2
        +\left(\frac{1}{\gamma^2}-1\right)
            \|\bar{\y}^k-\bar{\x}^k\|^2 \\
        &\qquad +2\left(\frac{1}{\gamma}-1\right)\bar{g}(\bar{\y}^k).
    \end{align*}
    We compute $V_k^\star-V_{k+1}^\star$ from this last expression,
    term by term. The update for $\bar{\z}^k$ gives
    \begin{equation}
        \label{eq:term_1}
        \|\bar{\z}^k\|^2-\|\bar{\z}^{k+1}\|^2
        =-2\theta_k\langle\bar{\z}^k,\bar{\y}^k-\bar{\x}^k\rangle
        -\theta_k^2\|\bar{\y}^k-\bar{\x}^k\|^2.
    \end{equation}
    \begin{equation}
        \label{eq:term_2}
            \langle\bar{\z}^k,\bar{\y}^k\rangle
            -\langle\bar{\z}^{k+1},\bar{\y}^{k+1}\rangle =-\theta_k\langle\bar{\y}^k-\bar{\x}^k,\bar{\y}^k\rangle
            -\langle\bar{\z}^{k+1},\bar{\y}^{k+1}-\bar{\y}^k\rangle.
    \end{equation}
    \begin{equation}
        \label{eq:term_3}
        -\|\bar{\y}^k\|^2+\|\bar{\y}^{k+1}\|^2
        =2\langle\bar{\y}^{k+1},\bar{\y}^{k+1}-\bar{\y}^k\rangle
        -\|\bar{\y}^{k+1}-\bar{\y}^k\|^2.
    \end{equation}
    Since $\bar{\v}^{k + 1} \in \partial\bar{g}(\bar{\y}^{k + 1})$, the definition of the generalized Bregman divergence gives
    \begin{equation}
        \label{eq:term_4}
        \bar{g}(\bar{\y}^k) - \bar{g}(\bar{\y}^{k + 1})
        = D_{\bar{g}}(\bar{\y}^k,\bar{\y}^{k + 1},\bar{\v}^{k + 1})
        - \langle \bar{\v}^{k + 1},\bar{\y}^{k + 1}-\bar{\y}^k\rangle.
    \end{equation}

    Using $\bar{\w}^{k+1}=\bar{\z}^{k+1}-\bar{\x}^{k+1}$,
    the definition of $\bar{\v}^{k+1}$ gives
    \[
        \bar{\v}^{k+1}
        =\bar{\y}^{k+1}-\bar{\z}^{k+1}
        -\left(1+\frac{1}{\gamma}\right)
        (\bar{\y}^{k+1}-\bar{\x}^{k+1}).
    \]

    Substituting this expression into \eqref{eq:term_4} gives
    \begin{equation}
    \begin{aligned}
        \label{eq:term_4_expanded}
        \bar{g}(\bar{\y}^k) - \bar{g}(\bar{\y}^{k + 1})
        &= D_{\bar{g}}(\bar{\y}^k,\bar{\y}^{k + 1},\bar{\v}^{k + 1}) - \langle \bar{\y}^{k + 1},\bar{\y}^{k + 1}-\bar{\y}^k\rangle + \langle \bar{\z}^{k + 1},\bar{\y}^{k + 1}-\bar{\y}^k\rangle
        \\ 
        &\qquad + \left(1 + \frac{1}{\gamma}\right)
        \langle \bar{\y}^{k + 1}-\bar{\x}^{k + 1},\bar{\y}^{k + 1}-\bar{\y}^k\rangle.
    \end{aligned}
    \end{equation}

    Combining these identities with their coefficients in $V_k^\star$, the terms involving $\langle\bar{\z}^{k + 1},\bar{\y}^{k + 1}-\bar{\y}^k\rangle$ and $\langle\bar{\y}^{k + 1},\bar{\y}^{k + 1}-\bar{\y}^k\rangle$ cancel. Adding $\left(\frac{1}{\gamma^2}-1\right)\left(\|\bar{\y}^k-\bar{\x}^k\|^2-\|\bar{\y}^{k+1}-\bar{\x}^{k+1}\|^2\right)$, we obtain
    \begin{equation}
    \label{eq:descent_intermediate}
    \begin{aligned}
        V_k^\star - V_{k + 1}^\star
        &= -2\theta_k\left\langle
        \bar{\z}^k + \left(\frac{1}{\gamma}-1\right)\bar{\y}^k,
        \bar{\y}^k-\bar{\x}^k\right\rangle
        -\theta_k^2\|\bar{\y}^k-\bar{\x}^k\|^2\\
        &\qquad +\left(\frac{1}{\gamma^2}-1\right)
        \bigg(\|\bar{\y}^k-\bar{\x}^k\|^2
        -\|\bar{\y}^{k + 1}-\bar{\x}^{k + 1}\|^2+2\langle\bar{\y}^{k + 1}-\bar{\x}^{k + 1},
        \bar{\y}^{k + 1}-\bar{\y}^k\rangle\bigg)\\
        &\qquad -\left(\frac{1}{\gamma}-1\right)
        \|\bar{\y}^{k + 1}-\bar{\y}^k\|^2+2\left(\frac{1}{\gamma}-1\right)
        D_{\bar{g}}(\bar{\y}^k,\bar{\y}^{k + 1},\bar{\v}^{k + 1}).
    \end{aligned}
    \end{equation}
    Using $\bar{\v}^k=(\bar{\x}^k-\bar{\y}^k)/\gamma-\bar{\w}^k$ and $\bar{\w}^k=\bar{\z}^k-\bar{\x}^k$, the first line of \eqref{eq:descent_intermediate} becomes
    \begin{equation}
    \label{eq:first_line_simplified}
    \begin{aligned}
        &-2\theta_k\left\langle\bar{\z}^k+\left(\frac{1}{\gamma}-1\right)\bar{\y}^k,\bar{\y}^k-\bar{\x}^k\right\rangle
        -\theta_k^2\|\bar{\y}^k-\bar{\x}^k\|^2 \\
        &\qquad =\theta_k(2-\theta_k)\|\bar{\y}^k-\bar{\x}^k\|^2
        +2\theta_k\langle\bar{\w}^k,\bar{\x}^k\rangle
        +2\theta_k\langle\bar{\v}^k,\bar{\y}^k\rangle.
    \end{aligned}
    \end{equation}
    To simplify the second line, use
    \[
        \bar{\y}^{k + 1}-\bar{\y}^k
        =(\bar{\y}^{k + 1}-\bar{\x}^{k + 1})-(\bar{\y}^k-\bar{\x}^k)
        +(\bar{\x}^{k + 1}-\bar{\x}^k).
    \]
    Expanding and collecting terms gives
    \begin{equation}
    \label{eq:residual_square_expansion}
    \begin{aligned}
        &\|\bar{\y}^k-\bar{\x}^k\|^2
        -\|\bar{\y}^{k + 1}-\bar{\x}^{k + 1}\|^2
        +2\langle\bar{\y}^{k + 1}-\bar{\x}^{k + 1},\bar{\y}^{k + 1}-\bar{\y}^k\rangle\\
        &\quad=2\langle\bar{\y}^k-\bar{\x}^k,\bar{\x}^{k + 1}-\bar{\x}^k\rangle
        -\|\bar{\x}^{k + 1}-\bar{\x}^k\|^2
        +\|\bar{\y}^{k + 1}-\bar{\y}^k\|^2.
    \end{aligned}
    \end{equation}
    Since $\bar{\z}^k=\bar{\x}^k+\bar{\w}^k$, the update \eqref{alg:step_3} implies
    \[
        \theta_k(\bar{\y}^k-\bar{\x}^k)
        = \bar{\z}^{k + 1} - \bar{\z}^k = (\bar{\x}^{k + 1}-\bar{\x}^k)+(\bar{\w}^{k + 1}-\bar{\w}^k).
    \]
    Dividing by $\theta_k$, taking the inner product with $2(\bar{\x}^{k+1}-\bar{\x}^k)$, and subtracting $\|\bar{\x}^{k+1}-\bar{\x}^k\|^2$ gives
    \begin{equation}
    \label{eq:successive_iterate_pairing}
        2\langle\bar{\y}^k-\bar{\x}^k,\bar{\x}^{k + 1}-\bar{\x}^k\rangle
        -\|\bar{\x}^{k + 1}-\bar{\x}^k\|^2 = \frac{2-\theta_k}{\theta_k}\|\bar{\x}^{k + 1}-\bar{\x}^k\|^2
        +\frac{2}{\theta_k}\langle\bar{\w}^{k + 1}-\bar{\w}^k,\bar{\x}^{k + 1}-\bar{\x}^k\rangle.
    \end{equation}
    We now return to the original coordinates. Applying the
    translation identity \eqref{eq:bregman_affine_translation}
    with $a=\x^\star$ and $b=\w^\star$, and using
    $\bar{\v}^{k+1}=\v^{k+1}+\w^\star$, gives
    \begin{equation}
        \label{eq:bregman_translation}
        D_{\bar{g}}(\bar{\y}^k,\bar{\y}^{k+1},\bar{\v}^{k+1})
        =D_g(\y^k,\y^{k+1},\v^{k+1}).
    \end{equation}
    The shifts cancel in $\bar{\y}^k-\bar{\x}^k$ and in
    differences between successive iterates, while the two
    pairings in \eqref{eq:first_line_simplified} satisfy
    \begin{align*}
        \langle\bar{\w}^k,\bar{\x}^k\rangle
        &=\langle\w^k-\w^\star,\x^k-\x^\star\rangle,\\
        \langle\bar{\v}^k,\bar{\y}^k\rangle
        &=\langle\v^k+\w^\star,\y^k-\x^\star\rangle.
    \end{align*}
    Substituting \eqref{eq:first_line_simplified},
    \eqref{eq:residual_square_expansion}, and
    \eqref{eq:successive_iterate_pairing} into
    \eqref{eq:descent_intermediate}, and using these translation
    identities, gives \eqref{eq:exact_descent_identity}. By convexity,
    the monotonicity inequalities established above, and the parameter
    bounds, every term on the right-hand side of
    \eqref{eq:exact_descent_identity} is nonnegative.
    Retaining the squared-residual term proves
    \eqref{eq:descent_bound}.
\end{proof}

The following lemma relates the normalized iteration with parameter $\gamma$ to a primal--dual iteration with parameter $\gamma^{-1}$ and rescaled relaxation parameters. This allows us to extend the analysis from $0<\gamma\le1$ to $\gamma>1$. Related primal--dual results include the classical self-duality identity for Douglas--Rachford splitting \cite[Corollary~4.3]{bauschke2012attouch} and iterate correspondences for compositions of relaxed resolvents \cite[Section~3]{alcantara2026douglas}.

\begin{lemma}[Primal--dual correspondence]
    \label{lem:edr_duality}
    Fix $\gamma>0$, let $(\theta_k)_{k\ge0}$ be a sequence of positive relaxation parameters, and let $(\x^k,\y^k,\z^k)_{k\ge0}$ be generated by \eqref{eq:edr_alpha_one}. Define
    \[
        \tilde{F}\coloneq F^{-1},\qquad \tilde{g}(u)\coloneq g^*(-u),\qquad
        \tilde{\gamma}\coloneq\frac{1}{\gamma},\qquad \tilde{\theta}_k\coloneq\gamma\theta_k\quad(k\ge0).
    \]

    Let $(\tilde{\x}^k,\tilde{\y}^k,\tilde{\z}^k)_{k\ge0}$ be generated by \eqref{eq:edr_alpha_one} with $(F,g,\gamma,\theta_k)$ replaced by $(\tilde{F},\tilde{g},\tilde{\gamma},\tilde{\theta}_k)$, starting from $\tilde{\z}^0=\z^0$. Then, for every $k\ge0$,
    \begin{equation}
        \label{eq:dual_iterate_correspondence}
        \tilde{\z}^k=\z^k,\qquad
        \tilde{\x}^k=\w^k,\qquad
        \tilde{\y}^k=-\v^k.
    \end{equation}
    Consequently,
    \[
        \tilde{\y}^k-\tilde{\x}^k
        =-\v^k-\w^k
        =\frac{\y^k-\x^k}{\gamma}.
    \]
    Moreover, for any $\x^\star,\w^\star\in\mathcal{H}$, the following equivalence holds.
    \begin{equation}
        \label{eq:dual_solution_pair_correspondence}
        \left\{\begin{aligned}
            \w^\star&\in F(\x^\star),\\
            -\w^\star&\in\partial g(\x^\star)
        \end{aligned}\right.
        \quad\Longleftrightarrow\quad
        \left\{\begin{aligned}
            \x^\star&\in \tilde{F}(\w^\star),\\
            -\x^\star&\in\partial \tilde{g}(\w^\star).
        \end{aligned}\right.
    \end{equation}
    In particular, the normalized inclusion has a solution if and only if the transformed inclusion has a solution.
\end{lemma}
\begin{proof}
    The subgradient inequality for $\tilde{g}(u) \coloneq g^*(-u)$ gives
    \[
        \partial\tilde{g}(u)=-\partial g^*(-u).
    \]
    Together with \eqref{eq:fenchel_subgradient}, this yields
    \begin{equation}
        \label{eq:dual_subgradient_correspondence}
        v\in\partial g(y)\quad\Longleftrightarrow\quad y\in\partial g^*(v)\quad\Longleftrightarrow\quad -y\in\partial\tilde{g}(-v).
    \end{equation}
    The solution-pair equivalence \eqref{eq:dual_solution_pair_correspondence} follows from $\tilde{F}=F^{-1}$ and the preceding subgradient equivalence. Thus interchanging the entries of a solution pair gives a solution pair for the transformed inclusion.

    We prove the iterate correspondence by induction. The initial states satisfy $\tilde{\z}^0=\z^0$. Suppose $\tilde{\z}^k=\z^k$.
    Since $\w^k=\z^k-\x^k\in F(\x^k)$, we have
    $\z^k-\w^k=\x^k\in\tilde{F}(\w^k)$.
    The resolvent characterization
    \eqref{eq:resolvent_characterization} therefore gives
    \[
        \tilde{\x}^k
        =J_{\tilde{F}}(\tilde{\z}^k)
        =\w^k.
    \]
    To show that $\tilde{\y}^k=-\v^k$, it suffices to verify that
    \[
        -\v^k
        =\operatorname{prox}_{\tilde{\gamma}\tilde{g}}
            \left(\tilde{\x}^k+
            \tilde{\gamma}(\tilde{\x}^k-\tilde{\z}^k)\right).
    \]
    Using $\tilde{\x}^k=\w^k$, $\tilde{\z}^k=\z^k$, and $\z^k=\x^k+\w^k$, the proximal argument reduces to $\w^k- \tilde{\gamma} \x^k$. The proximal optimality condition at $-\v^k$ is then verified by
    \[
        \frac{(\w^k-\tilde{\gamma}\x^k)-(-\v^k)}{\tilde{\gamma}}
        =\gamma(\w^k+\v^k)-\x^k
        =-\y^k\in\partial\tilde{g}(-\v^k).
    \]
    Here the last equality follows from the definition of $\v^k$,
    and the inclusion follows from
    \eqref{eq:dual_subgradient_correspondence}.
    Thus $\tilde{\y}^k=-\v^k$.
    Finally, since $\tilde{\theta}_k=\gamma\theta_k$ and $\gamma(\w^k+\v^k)=\x^k-\y^k$,
    \[
        \tilde{\z}^{k+1}=\tilde{\z}^k+\tilde{\theta}_k(\tilde{\y}^k-\tilde{\x}^k)=\z^k-\gamma\theta_k(\v^k+\w^k)=\z^k+\theta_k(\y^k-\x^k)=\z^{k+1}.
    \]
    This completes the induction and proves the iterate correspondence for every $k\ge0$.
\end{proof}

\begin{lemma}[Residual stability]
    \label{lem:edr_residual_stability}
    Fix $\gamma>0$ and let the iterates be generated by \eqref{eq:edr_alpha_one} with $0<\theta_k\le\min\{2,2/\gamma\}$. For every $0\le m\le N$,
    \begin{equation}
        \label{eq:edr_residual_stability}
        \|\y^N-\x^N\|^2\le\max\{\gamma,\gamma^{-1}\}\,\|\y^m-\x^m\|^2.
    \end{equation}
    The factor is sharp already on $\mathbb R^2$, for every constant relaxation parameter in the interior of the admissible interval.
\end{lemma}
\begin{proof}
    First suppose $0<\gamma\le1$. For $a\in\operatorname{dom}g^*$, define
    \begin{equation}\label{eq:residual_comparison}
        P_k(a)\coloneq\|\y^k-\x^k\|^2+
        \frac{\gamma(1-\gamma)}{(1+\gamma)^2}
        \left\{2D_{g^*}(a,\v^k,\y^k)+\|\y^k-\z^k-a\|^2\right\}.
    \end{equation}
    We prove the bound by showing that $\inf_{a\in\operatorname{dom}g^*}P_k(a)$ is nonincreasing and lies between $\|\y^k-\x^k\|^2$ and $\gamma^{-1}\|\y^k-\x^k\|^2$. Comparing its values at iterations $m$ and $N$ will then give the required factor $\gamma^{-1}$.

    The divergence is well defined because $\v^k\in\partial g(\y^k)$ implies $\y^k\in\partial g^*(\v^k)$ and, in particular, $\v^k\in\operatorname{dom}g^*$. Both terms inside the braces are nonnegative, so $P_k(a)\ge\|\y^k-\x^k\|^2$. For an upper bound on the infimum, take $a=\v^k$. The divergence vanishes, and the definition of $\v^k$ gives
    \[
        \y^k-\z^k-\v^k=\frac{1+\gamma}{\gamma}(\y^k-\x^k).
    \]
    Thus
    \[
        \|\y^k-\x^k\|^2
        \le\inf_{a\in\operatorname{dom}g^*}P_k(a)
        \le P_k(\v^k)=\gamma^{-1}\|\y^k-\x^k\|^2.
    \]
    It remains to prove that the infimum is nonincreasing.

    Fix $k\ge0$ and any $a\in\operatorname{dom}g^*$, and set
    \[
        b \coloneq a+\frac{\gamma\theta_k}{1+\gamma}(\v^k-a).
    \]
    Since $0<\gamma\theta_k/(1+\gamma)\le1$, this is a convex
    combination of $a$ and $\v^k$, so $b\in\operatorname{dom}g^*$.
    We estimate the changes in the three terms of $P_{k+1}(b)-P_k(a)$ separately and show that their sum is nonpositive.

    First, consider the squared-norm term inside the braces.
    The definitions of $\w^k$ and $\v^k$ give
    \[
        \theta_k(\y^k-\x^k)+(b-a)
        =\frac{\gamma\theta_k}{1+\gamma}(\y^k-\z^k-a).
    \]
    Substituting the update
    $\z^{k+1}=\z^k+\theta_k(\y^k-\x^k)$ therefore yields
    \begin{equation}\label{eq:residual_comparison_update}
        \y^{k+1}-\z^{k+1}-b
        =\left(1-\frac{\gamma\theta_k}{1+\gamma}\right)(\y^k-\z^k-a)
        +(\y^{k+1}-\y^k).
    \end{equation}
    Squaring this identity and subtracting $\| \y^k - \z^k - a \|^2$ from each side gives
    \begin{equation}\label{eq:other_squared_norm_change}
    \begin{aligned}
        &\|\y^{k+1}-\z^{k+1}-b\|^2-\|\y^k-\z^k-a\|^2\\
        &\quad=\left[\left(1-\frac{\gamma\theta_k}{1+\gamma}\right)^2-1\right]
            \|\y^k-\z^k-a\|^2\\
        &\qquad+2\left(1-\frac{\gamma\theta_k}{1+\gamma}\right)
            \langle\y^k-\z^k-a,\y^{k+1}-\y^k\rangle
            +\|\y^{k+1}-\y^k\|^2.
    \end{aligned}
    \end{equation}

    Next, we bound the change in the Bregman-divergence term. The definition of the Bregman divergence and Fenchel--Young equality \eqref{eq:fenchel_subgradient} give
    \[
        D_{g^*}(a,\v^k,\y^k)=g^*(a)+g(\y^k)-\langle a,\y^k\rangle.
    \]
    Subtracting this formula at $k$ from the same formula at
    $k+1$, and using the definition of $D_g$, gives
    \begin{equation}\label{eq:fixed_comparison_divergence}
        D_{g^*}(a,\v^{k+1},\y^{k+1})-D_{g^*}(a,\v^k,\y^k)
        =\langle\v^k-a,\y^{k+1}-\y^k\rangle+D_{g}(\y^{k+1},\y^k,\v^k).
    \end{equation}
    Since $b$ is a convex combination of $a$ and $\v^k$, convexity in the first argument gives
    \begin{equation}\label{eq:residual_divergence_convexity}
        D_{g^*}(b,\v^{k+1},\y^{k+1})
        \le \left(1-\frac{\gamma\theta_k}{1+\gamma}\right)D_{g^*}(a,\v^{k+1},\y^{k+1})+\frac{\gamma\theta_k}{1+\gamma}D_{g^*}(\v^k,\v^{k+1},\y^{k+1}).
    \end{equation}
    Since $D_{g^*}(a,\v^{k+1},\y^{k+1})\ge0$, we may increase its coefficient to $1$. For the second term, Fenchel--Young equality and the subgradient inequality at $\y^{k+1}$ give
    \begin{equation}\label{eq:successive_divergences}
        D_{g^*}(\v^k,\v^{k+1},\y^{k+1})=D_{g}(\y^{k+1},\y^k,\v^k)
        \le\langle\v^{k+1}-\v^k,\y^{k+1}-\y^k\rangle.
    \end{equation}
    Subtract $D_{g^*}(a,\v^k,\y^k)$ from both sides of \eqref{eq:residual_divergence_convexity}, with its first coefficient replaced by $1$. Substituting \eqref{eq:fixed_comparison_divergence} and the equality in \eqref{eq:successive_divergences} then gives
    \[
    \begin{aligned}
        &D_{g^*}(b,\v^{k+1},\y^{k+1})-D_{g^*}(a,\v^k,\y^k)\\
        &\qquad\le\langle\v^k-a,\y^{k+1}-\y^k\rangle
        +\left(1+\frac{\gamma\theta_k}{1+\gamma}\right)
        D_g(\y^{k+1},\y^k,\v^k).
    \end{aligned}
    \]
    Applying the inequality in \eqref{eq:successive_divergences} to the final term yields
    \begin{equation}\label{eq:residual_divergence_change}
    \begin{aligned}
        &D_{g^*}(b,\v^{k+1},\y^{k+1})-D_{g^*}(a,\v^k,\y^k)\\
        &\qquad\le\langle\v^k-a,\y^{k+1}-\y^k\rangle+
        \left(1+\frac{\gamma\theta_k}{1+\gamma}\right)
        \langle\v^{k+1}-\v^k,\y^{k+1}-\y^k\rangle.
    \end{aligned}
    \end{equation}

    Finally, we bound the change in the squared residual. We claim that
    \begin{equation}\label{eq:squared_residual_inequality}
    \begin{aligned}
        &\frac{\theta_k}{2}\bigl(\|\y^{k+1}-\x^{k+1}\|^2-\|\y^k-\x^k\|^2\bigr)\\
        &\qquad\le\left\langle\y^{k+1}-\y^k,\,
        (\y^{k+1}-\x^{k+1})-(\y^k-\x^k)+\w^{k+1}-\w^k\right\rangle.
    \end{aligned}
    \end{equation}
    To verify this inequality, expand the squared-norm difference
    on the left and the inner product on the right. Using the
    update identity
    $\theta_k(\y^k-\x^k)=(\x^{k+1}-\x^k)+(\w^{k+1}-\w^k)$,
    the right-hand side minus the left-hand side is
    \[
        \left(1-\frac{\theta_k}{2}\right)
        \bigl\|(\y^{k+1}-\x^{k+1})-(\y^k-\x^k)\bigr\|^2
        +\langle\x^{k+1}-\x^k,\w^{k+1}-\w^k\rangle.
    \]
    The first term is nonnegative because $\theta_k\le2$, and the second by monotonicity of $F$.
    Rewriting the right-hand side of \eqref{eq:squared_residual_inequality} using $\v^k=(\x^k-\y^k)/\gamma-\w^k$ and the $\z$-update \eqref{alg:step_3}, then dividing by $\theta_k/2>0$, gives
    \begin{equation}\label{eq:residual_pairing}
    \begin{aligned}
        &\|\y^{k+1}-\x^{k+1}\|^2-\|\y^k-\x^k\|^2\\
        &\quad\le\frac{2(1-\gamma)}{\theta_k(1+\gamma)}
        \left\langle\y^{k+1}-\y^k,\,\theta_k(\y^k-\x^k)-(\y^{k+1}-\y^k)\right\rangle\\
        &\qquad-\frac{4\gamma}{\theta_k(1+\gamma)}
        \langle\v^{k+1}-\v^k,\y^{k+1}-\y^k\rangle.
    \end{aligned}
    \end{equation}

    We combine \eqref{eq:residual_pairing}, \eqref{eq:residual_divergence_change}, and \eqref{eq:other_squared_norm_change} with their coefficients in $P_k$. Using the identity
    \[
        (\y^k-\x^k)+\frac{\gamma}{1+\gamma}(\v^k-a)
        =\frac{\gamma}{1+\gamma}(\y^k-\z^k-a)
    \]
    to collect the mixed terms and complete the square gives
    \begin{equation}\label{eq:combined_residual_comparison}
    \begin{aligned}
        P_{k+1}(b)-P_k(a)
        &\le\left[-\frac{4\gamma}{\theta_k(1+\gamma)}
        +\frac{2\gamma(1-\gamma)}{(1+\gamma)^2}
        \left(1+\frac{\gamma\theta_k}{1+\gamma}\right)\right]\langle\v^{k+1}-\v^k,\y^{k+1}-\y^k\rangle\\
        &\qquad-\frac{1-\gamma}{1+\gamma}
        \left(\frac{2}{\theta_k}-\frac{\gamma}{1+\gamma}\right)\left\|\y^{k+1}-\y^k-\frac{\gamma\theta_k}{1+\gamma}(\y^k-\z^k-a)\right\|^2\le0.
    \end{aligned}
    \end{equation}
    Since $0<\gamma\le1$ and $0<\theta_k\le2$, the coefficient of $\langle\v^{k+1}-\v^k,\y^{k+1}-\y^k\rangle$ in \eqref{eq:combined_residual_comparison} is at most
    \[
        -\frac{2\gamma}{1+\gamma}
        +\frac{2\gamma(1-\gamma)(1+3\gamma)}{(1+\gamma)^3}
        =-\frac{8\gamma^3}{(1+\gamma)^3}\le0.
    \]
    Monotonicity of $\partial g$ makes this inner product nonnegative, so the first term is nonpositive. The second is nonpositive because $1-\gamma\ge0$ and $2/\theta_k-\gamma/(1+\gamma)\ge1/(1+\gamma)>0$.

    Since $b$ is admissible, the next infimum is at most $P_{k+1}(b)$, and hence at most $P_k(a)$. This holds for every $a\in\operatorname{dom}g^*$, so taking the infimum over $a$ gives
    \[
        \inf_{a\in\operatorname{dom}g^*}P_{k+1}(a)
        \le\inf_{a\in\operatorname{dom}g^*}P_k(a).
    \]
    The lower bound at $N$ and upper bound at $m$ therefore give
    \begin{equation}
        \label{eq:residual_stability_small_gamma}
        \|\y^N-\x^N\|^2\le\gamma^{-1}\|\y^m-\x^m\|^2.
    \end{equation}
    This includes $\gamma=1$, when $P_k(a)$ is the squared residual for every $a\in\operatorname{dom}g^*$.

    For $\gamma>1$, Lemma~\ref{lem:edr_duality} gives
    \[
        0<\tilde{\gamma}=\gamma^{-1}<1,\qquad
        0<\tilde{\theta}_k=\gamma\theta_k\le2,\qquad
        \tilde{\y}^k-\tilde{\x}^k=\frac{\y^k-\x^k}{\gamma}.
    \]
    Thus \eqref{eq:residual_stability_small_gamma} applies to the tilded iterates and gives
    \[
        \|\tilde{\y}^N-\tilde{\x}^N\|^2
        \le\tilde{\gamma}^{-1}\|\tilde{\y}^m-\tilde{\x}^m\|^2.
    \]
    Substituting $\tilde{\gamma}^{-1}=\gamma$ and
    $\tilde{\y}^k-\tilde{\x}^k=(\y^k-\x^k)/\gamma$
    for $k=m,N$, and multiplying by $\gamma^2$, gives
    \begin{equation}
        \label{eq:residual_stability_large_gamma}
        \|\y^N-\x^N\|^2\le\gamma\|\y^m-\x^m\|^2.
    \end{equation}
    Combining \eqref{eq:residual_stability_small_gamma} and \eqref{eq:residual_stability_large_gamma} proves \eqref{eq:edr_residual_stability}.

    To see optimality, fix $0<\theta<\min\{2,2/\gamma\}$ and take $F=\partial\iota_U$ and $g=\iota_V$, where $U$ is the horizontal axis in $\mathbb R^2$ and $V$ is the line making angle $\varepsilon\in(0,\pi/2)$ with $U$. Both steps are orthogonal projections. Substituting them into \eqref{eq:edr_alpha_one} gives $\y^k-\x^k=R_\varepsilon\z^k$, where
    \[
        R_\varepsilon=\sin\varepsilon
        \begin{pmatrix}-\sin\varepsilon&-\gamma\cos\varepsilon\\
                         \cos\varepsilon&-\gamma\sin\varepsilon\end{pmatrix}.
    \]
    With $\theta_k=\theta$, the update \eqref{alg:step_3} gives
    \[
        \y^{k+1}-\x^{k+1}
        =R_\varepsilon\left(\z^k+\theta(\y^k-\x^k)\right)
        =(I+\theta R_\varepsilon)(\y^k-\x^k).
    \]
    Iterating this identity gives
    \begin{equation}
        \label{eq:sharpness_residual_evolution}
        \y^N-\x^N=(I+\theta R_\varepsilon)^N(\y^0-\x^0).
    \end{equation}
    Since $\sin\varepsilon/\varepsilon\to1$ as $\varepsilon\downarrow0$, we have
    \[
        \frac{R_\varepsilon}{\varepsilon}\to
        B\coloneq\begin{pmatrix}0&-\gamma\\1&0\end{pmatrix}.
    \]
    To obtain a limiting matrix power, choose
    $N=\lfloor\pi/(2\theta\sqrt{\gamma}\,\varepsilon)\rfloor$.
    Then $N\theta\varepsilon\to\pi/(2\sqrt{\gamma})$, and hence
    $N\theta R_\varepsilon\to\pi B/(2\sqrt{\gamma})$.
    Since $N\to\infty$ as $\varepsilon\downarrow0$,
    the matrix exponential limit gives
    \begin{equation}
        \label{eq:sharpness_matrix_power_limit}
        (I+\theta R_\varepsilon)^N\to
        \exp\!\left(\frac{\pi B}{2\sqrt\gamma}\right)=\frac{B}{\sqrt\gamma}.
    \end{equation}
    The last equality follows from $B^2=-\gamma I$.
    Since $\det R_\varepsilon=\gamma\sin^2\varepsilon>0$,
    $R_\varepsilon$ is invertible. We can therefore choose $\z^0$
    so that the initial residual $\y^0-\x^0=R_\varepsilon\z^0$
    is either coordinate unit vector.
    Combining \eqref{eq:sharpness_residual_evolution} and
    \eqref{eq:sharpness_matrix_power_limit} gives, as
    $\varepsilon\downarrow0$,
    \[
        \y^N-\x^N\to\frac{B}{\sqrt{\gamma}}(\y^0-\x^0).
    \]
    Taking squared norms for the two choices of initial residual yields
    \begin{equation}
        \label{eq:sharpness_residual_ratios}
        \frac{\|\y^N-\x^N\|^2}{\|\y^0-\x^0\|^2}
        \to
        \begin{cases}
            \gamma^{-1}, & \y^0-\x^0=(1,0)^\top,\\
            \gamma, & \y^0-\x^0=(0,1)^\top.
        \end{cases}
    \end{equation}
    Any constant replacing $\max\{\gamma,\gamma^{-1}\}$ in
    \eqref{eq:edr_residual_stability} must bound the left-hand side
    of \eqref{eq:sharpness_residual_ratios} for both initializations,
    for every $\varepsilon>0$ and its corresponding $N$.
    Passing to the limit forces that constant to be at least
    $\max\{\gamma,\gamma^{-1}\}$, proving sharpness.
\end{proof}

\subsection{Proof of Theorem~\ref{thm:edr_convergence}}

We now combine the three lemmas to prove
Theorem~\ref{thm:edr_convergence}.
For $0<\gamma\le1$, we establish boundedness of
$(\x^k)$, $(\y^k)$, $(\z^k)$, and $(\w^k)$ and derive
the residual rate \eqref{eq:edr_residual_rate}.
We then show that every weak cluster point of $(\x^k,\w^k)$
is a solution pair and that there is only one such cluster point.
Lemma~\ref{lem:edr_duality} extends these conclusions to $\gamma>1$.

\begin{proof}
    Returning to the notation of Theorem~\ref{thm:edr_convergence}, we set $\gamma=\beta/\alpha$ and apply the three lemmas with $(\alpha F,\alpha g)$ in place of $(F,g)$. Since $\partial(\alpha g)=\alpha\partial g$, we have
    \[
        \operatorname{zer}(\alpha F+\partial(\alpha g))
        =\operatorname{zer}(F+\partial g)\ne\varnothing.
    \]
    For these scaled operators, $\w^k=\z^k-\x^k\in\alpha F(\x^k)$ and $\v^k\in\alpha\partial g(\y^k)$. The relaxation bound becomes $\vartheta=\min\{2,2/\gamma\}$.

    First suppose $0<\gamma\le1$, so $\vartheta=2$. Choose $\x^\star,\w^\star$ satisfying
    \[
        \w^\star\in\alpha F(\x^\star),\qquad
        -\w^\star\in\alpha\partial g(\x^\star).
    \]
    By Lemma~\ref{lem:edr_descent}, the Lyapunov function $V_k^\star$ centered at this solution pair is nonnegative and nonincreasing, so it converges.

    To establish boundedness of $(\x^k)$ and $(\w^k)$, we rewrite
    the two norm terms in \eqref{eq:lyapunov} containing
    $\y^k-\x^k$ as
    \[
    \begin{aligned}
        &\left\|\w^k-\w^\star+
        \left(\frac1\gamma-1\right)(\y^k-\x^k)\right\|^2
        +\left(\frac1\gamma-1\right)\|\y^k-\x^k\|^2\\
        &\qquad=\gamma\|\w^k-\w^\star\|^2
        +(1-\gamma)\left\|\w^k-\w^\star+
        \frac{\y^k-\x^k}{\gamma}\right\|^2.
    \end{aligned}
    \]
    Substituting this identity into \eqref{eq:lyapunov} and discarding
    the remaining nonnegative terms gives
    \begin{equation}
        \label{eq:convergence_boundedness}
        \frac1\gamma\|\x^k-\x^\star\|^2+\gamma\|\w^k-\w^\star\|^2
        \le V_k^\star\le V_0^\star.
    \end{equation}
    Thus $(\x^k)$ and $(\w^k)$ are bounded, and so is
    $\z^k=\x^k+\w^k$. Nonexpansiveness of the proximal operator yields
    \[
        \|\y^k-\x^\star\|
        =\|\operatorname{prox}_{\beta g}(\x^k-\gamma\w^k)
          -\operatorname{prox}_{\beta g}(\x^\star-\gamma\w^\star)\|
        \le\|(\x^k-\x^\star)-\gamma(\w^k-\w^\star)\|,
    \]
    proving boundedness of $(\y^k)$.

    To prove \eqref{eq:edr_residual_rate}, fix $m\ge0$. For every $N>m$, apply Lemma~\ref{lem:edr_residual_stability} from each $k=m,\ldots,N-1$ to $N$ and sum with weights $\theta_k(2-\theta_k)$. Summing \eqref{eq:descent_bound} over the same indices gives
    \begin{equation}
    \label{eq:residual_tail_bound}
        \left(\sum_{k=m}^{N-1}\theta_k(2-\theta_k)\right)\|\y^N-\x^N\|^2\le \frac1\gamma\sum_{k=m}^{N-1}\theta_k(2-\theta_k)\|\y^k-\x^k\|^2\le \frac{V_m^\star-V_N^\star}{\gamma}.
    \end{equation}
    Removing the first $m$ terms from the divergent sum
    in \eqref{eq:edr-parameters} changes it by a finite amount
    independent of $N$, so
    \[
        \frac{\sum_{k=m}^{N-1}\theta_k(2-\theta_k)}
             {\sum_{k=0}^{N-1}\theta_k(2-\theta_k)}
        \to1
        \qquad\text{as }N\to\infty.
    \]
    Multiplying \eqref{eq:residual_tail_bound} by the reciprocal
    of this ratio and taking the limit superior in $N$ gives
    \[
        0\le
        \limsup_{N\to\infty}
        \left(\sum_{k=0}^{N-1}\theta_k(2-\theta_k)\right)
        \|\y^N-\x^N\|^2
        \le
        \frac{V_m^\star-\lim_{N\to\infty}V_N^\star}{\gamma}.
    \]
    The limit superior is independent of $m$.
    Since $(V_m^\star)$ converges, letting $m\to\infty$ makes the
    right-hand side tend to zero. The limit superior is therefore zero, which proves
    \eqref{eq:edr_residual_rate}, with $\vartheta=2$.
    Since the sum of weights diverges, this also implies
    $\y^N-\x^N\to0$.

    We next show that every weak cluster point of $(\x^k,\w^k)$
    is a solution pair.
    By \eqref{eq:convergence_boundedness}, $(\x^k,\w^k)$ is bounded
    and therefore has a weak cluster point.
    Let $(\x^\infty,\w^\infty)$ be any such point and choose
    a subsequence such that
    \[
        (\x^{k_j},\w^{k_j})\rightharpoonup(\x^\infty,\w^\infty).
    \]
    Since $\y^k-\x^k\to0$, we obtain
    \[
        \y^{k_j}\rightharpoonup\x^\infty,\qquad
        \v^{k_j}\rightharpoonup-\w^\infty,\qquad
        \w^k+\v^k\to0.
    \]
    By boundedness and the strong convergences established above,
    \begin{equation}
        \label{eq:graph_pairing_limit}
        \langle \x^k,\w^k\rangle+\langle \y^k,\v^k\rangle=\langle \x^k-\y^k,\w^k\rangle+\langle \y^k,\w^k+\v^k\rangle\to0.
    \end{equation}
    Take arbitrary $u\in\alpha F(a)$ and $v\in\alpha\partial g(b)$. Monotonicity gives
    \begin{equation}
        \label{eq:cluster_monotonicity_test}
        \langle \x^{k_j}-a,\w^{k_j}-u\rangle+\langle \y^{k_j}-b,\v^{k_j}-v\rangle\ge0.
    \end{equation}
    Expanding \eqref{eq:cluster_monotonicity_test}, the terms
    $\langle\x^{k_j},\w^{k_j}\rangle$ and
    $\langle\y^{k_j},\v^{k_j}\rangle$ have sum tending to zero
    by \eqref{eq:graph_pairing_limit}.
    Every remaining inner product contains a fixed vector among
    $a$, $b$, $u$, and $v$, so weak convergence allows passage
    to the limit. Regrouping gives
    \begin{equation}
        \label{eq:limiting_monotonicity_test}
        \langle \x^\infty-a,\w^\infty-u\rangle+\langle \x^\infty-b,-\w^\infty-v\rangle\ge0.
    \end{equation}
    Thus $\bigl((\x^\infty,\x^\infty),(\w^\infty,-\w^\infty)\bigr)$ is monotonically related to the graph of $(\alpha F)\times\partial(\alpha g)$. Since this product operator is maximally monotone, \eqref{eq:maximal_graph_characterization} gives
    \begin{equation}
        \label{eq:cluster_solution_pair}
        \w^\infty\in\alpha F(\x^\infty),\qquad -\w^\infty\in\alpha\partial g(\x^\infty).
    \end{equation}

    To prove that $(\x^k,\w^k)$ has a unique weak cluster point, let $(\x_1^\infty,\w_1^\infty)$ and $(\x_2^\infty,\w_2^\infty)$ be two weak cluster points. Both are solution pairs by \eqref{eq:cluster_solution_pair}. For $i=1,2$, let $V_k^{(i)}$ denote the Lyapunov function centered at $(\x_i^\infty,\w_i^\infty)$. Both sequences converge by Lemma~\ref{lem:edr_descent}.

    To compute $V_k^{(1)}-V_k^{(2)}$, expand \eqref{eq:lyapunov} at each center $(\x_i^\infty,\w_i^\infty)$ and use
    \[
        D_{\alpha g}(\y^k,\x_i^\infty,-\w_i^\infty)
        =\alpha g(\y^k) - \alpha g(\x_i^\infty)+\langle \w_i^\infty,\y^k\rangle-\langle \w_i^\infty,\x_i^\infty\rangle.
    \]
    All terms independent of the center cancel in the difference.
    The terms pairing $\y^k$ with $\w_i^\infty$ also cancel.
    Using $\z^k=\x^k+\w^k$, the remaining terms involving both
    the center and an iterate combine into
    \[
        -\frac{2}{\gamma}\langle\z^k,\x_i^\infty\rangle
        -2\langle\z^k,\w_i^\infty\rangle
        =-2\left\langle\z^k,
            \frac{\x_i^\infty}{\gamma}+\w_i^\infty\right\rangle.
    \]
    Collecting the terms depending only on the center and
    subtracting the expansions for $i=1$ and $i=2$ yields
    \begin{equation}
    \label{eq:lyapunov_center_difference}
    \begin{aligned}
        V_k^{(1)}-V_k^{(2)}
        ={}&\frac{\|\x_1^\infty\|^2-\|\x_2^\infty\|^2}{\gamma}+\|\w_1^\infty\|^2-\|\w_2^\infty\|^2\\
        &+2\bigl(\langle \x_1^\infty,\w_1^\infty\rangle-\langle \x_2^\infty,\w_2^\infty\rangle\bigr)
        -2\alpha\left(\frac1\gamma-1\right)\bigl(g(\x_1^\infty)-g(\x_2^\infty)\bigr)\\
        &-2\left\langle\z^k,\frac{\x_1^\infty-\x_2^\infty}{\gamma}+\w_1^\infty-\w_2^\infty\right\rangle.
    \end{aligned}
    \end{equation}

    Only the last inner product on the right-hand side of \eqref{eq:lyapunov_center_difference} depends on $k$. Since both Lyapunov sequences converge, this inner product has a limit. Along a subsequence converging weakly to $(\x_i^\infty,\w_i^\infty)$, we have $\z^k\rightharpoonup \x_i^\infty+\w_i^\infty$. Equating the two subsequential limits yields
    \begin{equation}
        \label{eq:cluster_point_identity}
        \left\langle (\x_1^\infty-\x_2^\infty)+(\w_1^\infty-\w_2^\infty),\frac{\x_1^\infty-\x_2^\infty}{\gamma}+(\w_1^\infty-\w_2^\infty)\right\rangle=0.
    \end{equation}
    Since $\w_i^\infty\in\alpha F(\x_i^\infty)$ and $-\w_i^\infty\in\alpha\partial g(\x_i^\infty)$, monotonicity of the two operators gives both $\langle \x_1^\infty-\x_2^\infty,\w_1^\infty-\w_2^\infty\rangle\ge0$ and $\langle \x_1^\infty-\x_2^\infty,\w_1^\infty-\w_2^\infty\rangle\le0$. The cross terms in \eqref{eq:cluster_point_identity} therefore vanish, leaving
    \[
        \frac1\gamma\|\x_1^\infty-\x_2^\infty\|^2+\|\w_1^\infty-\w_2^\infty\|^2=0.
    \]
    Thus $(\x^k,\w^k)$ has a unique weak cluster point.

    Since the sequence is bounded, the criterion in
    Section~\ref{sec:weak_convergence} implies convergence
    of the whole sequence:
    \[
        (\x^k,\w^k)\rightharpoonup(\x^\infty,\w^\infty).
    \]
    Since $\z^k=\x^k+\w^k$, we also obtain
    \[
        \z^k\rightharpoonup\z^\infty\coloneq\x^\infty+\w^\infty.
    \]
    Since $\y^k-\x^k\to0$, we also have
    $\y^k\rightharpoonup\x^\infty$.

    Since $\z^\infty-\x^\infty=\w^\infty$ and
    $\beta=\alpha\gamma$, the solution-pair inclusions
    \eqref{eq:cluster_solution_pair} give the limit identities
    \eqref{eq:edr_limit_identities} by the resolvent
    characterization \eqref{eq:resolvent_characterization}
    and the proximal optimality condition
    \eqref{eq:proximal_optimality}.
    Adding the inclusions and dividing by $\alpha$ gives
    $0\in F(\x^\infty)+\partial g(\x^\infty)$.

    Now suppose $\gamma>1$. Apply Lemma~\ref{lem:edr_duality} to the scaled operators $\alpha F,\alpha g$, defining
    \[
        \tilde{F}\coloneq(\alpha F)^{-1},\qquad
        \tilde{g}(u)\coloneq(\alpha g)^*(-u),\qquad
        \tilde{\gamma}\coloneq\frac1\gamma,\qquad
        \tilde{\theta}_k\coloneq\gamma\theta_k\quad(k\ge0).
    \]
    The transformed operators satisfy the same assumptions and admit a solution pair. Moreover, the parameter assumptions give
    \[
        0<\tilde{\gamma}<1,\qquad 0<\tilde{\theta}_k\le2\quad(k\ge0),
    \]
    and, since $\vartheta=2/\gamma$,
    \begin{equation}
        \label{eq:dual_damping_sum}
        \sum_{k=0}^{\infty}\tilde{\theta}_k(2-\tilde{\theta}_k)
        =\gamma^2\sum_{k=0}^{\infty}\theta_k(\vartheta-\theta_k)=\infty.
    \end{equation}
    Applying the case $0<\gamma\le1$ to the transformed iteration,
    initialized at $\tilde{\z}^0=\z^0$, gives weak convergence
    and the residual rate \eqref{eq:edr_residual_rate}
    for the transformed parameters.
    In particular,
    \[
        (\tilde{\x}^k,\tilde{\z}^k)\rightharpoonup(\tilde{\x}^\infty,\tilde{\z}^\infty),\qquad \tilde{\y}^k-\tilde{\x}^k\to0.
    \]
    By the iterate correspondence in Lemma~\ref{lem:edr_duality},
    \begin{equation}
        \label{eq:dual_iterate_recovery}
        \z^k=\tilde{\z}^k,\qquad \x^k=\tilde{\z}^k-\tilde{\x}^k,\qquad \y^k-\x^k=\gamma(\tilde{\y}^k-\tilde{\x}^k).
    \end{equation}
    The finite-sum version of \eqref{eq:dual_damping_sum} and
    the residual identity in \eqref{eq:dual_iterate_recovery} give
    \[
        \left(\sum_{k=0}^{N-1}\theta_k(\vartheta-\theta_k)\right)
        \|\y^N-\x^N\|^2
        =
        \left(\sum_{k=0}^{N-1}
        \tilde{\theta}_k(2-\tilde{\theta}_k)\right)
        \|\tilde{\y}^N-\tilde{\x}^N\|^2
        \to0.
    \]
    The first two identities in \eqref{eq:dual_iterate_recovery} also give
    \[
        (\x^k,\z^k)\rightharpoonup(\x^\infty,\z^\infty),
    \]
    where $\x^\infty\coloneq\tilde{\z}^\infty-\tilde{\x}^\infty$ and $\z^\infty\coloneq\tilde{\z}^\infty$.
    Since $\y^k-\x^k\to0$ by \eqref{eq:dual_iterate_recovery},
    we also have $\y^k\rightharpoonup\x^\infty$.
    The solution-pair correspondence
    \eqref{eq:dual_solution_pair_correspondence} gives
    \[
        \z^\infty-\x^\infty\in\alpha F(\x^\infty),
        \qquad
        -(\z^\infty-\x^\infty)\in\alpha\partial g(\x^\infty).
    \]
    As above, these inclusions imply the limit identities
    \eqref{eq:edr_limit_identities} and
    $0\in F(\x^\infty)+\partial g(\x^\infty)$.
    This completes the proof.
\end{proof}

\subsection{Necessity of the relaxation conditions}

\begin{proposition}
    \label{prop:edr_parameter_sharpness}
    Fix $\alpha,\beta>0$ and let
    $\vartheta=\min\{2,2\alpha/\beta\}$.
    Consider either of the following relaxation sequences:
    \begin{enumerate}[label=\roman*.]
        \item \label{item:sharpness_constant} constant relaxation $\theta_k=\theta$ with
        $\theta\notin(0,\vartheta)$;
        \item \label{item:sharpness_sequence} a prescribed sequence satisfying
        $0<\theta_k\le\vartheta$ for every $k\ge0$ and
        \begin{equation}
            \label{eq:sharpness_finite_damping}
            \sum_{k=0}^{\infty}
            \theta_k(\vartheta-\theta_k)<\infty.
        \end{equation}
    \end{enumerate}
    In either case, there exist a maximally monotone operator
    $F:\mathbb R^2\rightrightarrows\mathbb R^2$ and a proper,
    lower semicontinuous, convex function
    $g:\mathbb R^2\to(-\infty,+\infty]$ such that
    \[
        \operatorname{zer}(F+\partial g)=\{0\},
    \]
    but, for every $\z^0\ne0$, the iterates generated by
    \eqref{eq:extended_douglas_rachford} satisfy
    \begin{equation}
        \label{eq:sharpness_positive_distance}
        \inf_{k\ge0}\|\x^k\|>0.
    \end{equation}
\end{proposition}

\begin{proof}
Fix $\alpha,\beta>0$, write $\gamma=\beta/\alpha$, and let $R$
be rotation by $\pi/2$ on $\mathbb{R}^2$.
Take $F=R/\alpha$ and choose
\[
    g=\begin{cases}
        \iota_{\{0\}}, & \gamma\le1,\\
        0, & \gamma>1.
    \end{cases}
\]
The rotation $R$ is invertible and satisfies
$\langle x,Rx\rangle=0$ for every $x$.
Thus $F$ is monotone, and invertibility of
$I+\alpha F=I+R$ gives maximal monotonicity by Minty's theorem
\cite[Theorem~21.1]{bauschke2017convex}.
When $\gamma\le1$, every solution must lie in
$\operatorname{dom}\partial g=\{0\}$.
When $\gamma>1$, we have $\partial g=0$, so the inclusion
reduces to $Rx=0$.
Thus zero is the unique solution in both cases.

The resolvent step gives $\z^k=\x^k+R\x^k$. These two summands are orthogonal and have equal norm, so $\|\z^k\|^2=2\|\x^k\|^2$. If $\gamma\le1$, the proximal step gives $\y^k=0$. If $\gamma>1$, it gives $\y^k=\x^k-\gamma R\x^k$. Hence
\[
    \z^{k+1}=\begin{cases}
        (1-\theta_k)\x^k+R\x^k, & \gamma\le1,\\
        \x^k+(1-\gamma\theta_k)R\x^k, & \gamma>1.
    \end{cases}
\]
Writing $t_k=\max\{1,\gamma\}\theta_k$, we obtain from orthogonality in both cases
\begin{equation}
    \label{eq:sharpness_norm_recursion}
    \|\z^{k+1}\|^2
    =\frac{1+(1-t_k)^2}{2}\|\z^k\|^2
    =\left(1-\frac{1}{2} t_k(2-t_k)\right)\|\z^k\|^2.
\end{equation}
Since $\vartheta=2/\max\{1,\gamma\}$, constant relaxation
$\theta\notin(0,\vartheta)$ gives $t_k\notin(0,2)$.
The multiplier in \eqref{eq:sharpness_norm_recursion} is then
at least one, so for every $\z^0\ne0$,
\[
    \|\x^k\|^2=\frac{\|\z^k\|^2}{2}
    \ge\frac{\|\z^0\|^2}{2}>0.
\]
This proves \eqref{eq:sharpness_positive_distance}
in case~\ref{item:sharpness_constant}.

We next consider case~\ref{item:sharpness_sequence}.
Fix any prescribed sequence with $0<\theta_k\le\vartheta$.
Then $0<t_k\le2$ and
\begin{equation}
    \label{eq:sharpness_damping_scaling}
    t_k(2-t_k)=\max\{1,\gamma\}^2\theta_k(\vartheta-\theta_k).
\end{equation}
By \eqref{eq:sharpness_finite_damping} and
\eqref{eq:sharpness_damping_scaling}, we have
$\sum_k t_k(2-t_k)<\infty$.
Iterating \eqref{eq:sharpness_norm_recursion} gives
\begin{equation}
    \label{eq:sharpness_product}
    \|\z^n\|^2=\|\z^0\|^2\prod_{k=0}^{n-1}\left(1-\frac{1}{2} t_k(2-t_k)\right).
\end{equation}
Since $0\le\frac{1}{2} t_k(2-t_k)\le\frac{1}{2}$, we have
\[
    0\le-\log\left(1-\frac{1}{2} t_k(2-t_k)\right)
    \le t_k(2-t_k).
\]
Since $\sum_k t_k(2-t_k)<\infty$, the logarithm bound
shows that the logarithm of the product in
\eqref{eq:sharpness_product} converges to a finite real number.
The product therefore has a positive limit.
For every $\z^0\ne0$,
$\|\x^n\|^2=\|\z^n\|^2/2$ stays bounded away from zero,
proving \eqref{eq:sharpness_positive_distance}
in case~\ref{item:sharpness_sequence}.
\end{proof}

Together, Theorem~\ref{thm:edr_convergence} and
Proposition~\ref{prop:edr_parameter_sharpness} give necessary
and sufficient conditions for a universal convergence guarantee,
both for constant relaxation and for prescribed relaxation
sequences within the pointwise bounds.

\section{Conclusion}

We have established weak convergence of extended Douglas--Rachford splitting with independently chosen resolvent stepsizes and variable relaxation when the first operator is maximally monotone and the second is a convex subdifferential. Assuming the inclusion has a solution, our relaxation criterion extends the classical Douglas--Rachford condition and permits adaptive choices without requiring the relaxation parameters to converge. For prescribed sequences within the pointwise bounds, the divergence condition is also necessary for a universal convergence guarantee. The constant-relaxation specialization resolves the conjecture of Nilsson, Åkerman, and Giselsson.

The analysis combines a Lyapunov function whose definition is independent of the relaxation parameter with a solution-independent functional controlling residual growth. The resulting sharp stability bound quantifies how much the residual norm can increase when the resolvent stepsizes differ, and remains valid even when the inclusion has no solution. Under the convergence hypotheses, this control yields a last-iterate residual rate despite possible increases in the residual norm, recovering $o(N^{-1/2})$ when the relaxation parameters remain uniformly separated from both endpoints.

\printbibliography

\end{document}